\documentclass[12pt,oneside]{amsart}
\usepackage{cite}
\usepackage{calc}
\usepackage{setspace}
\usepackage[letterpaper,left=2cm,right=2cm,top=3cm,bottom=3cm]{geometry}

\def\nlss{nonlinear Schr\"odinger equations}

\def\U{\mathcal{U}}
\def\N{\mathbb{N}}
\def\C{\mathbb{C}}
\def\R{\mathbb{R}}
\def\B{\mathcal{B}}

\def\half{\frac{1}{2}}
\renewcommand{\epsilon}{\varepsilon}

\newtheorem{theorem}{Theorem}

\newtheorem{definition}[theorem]{Definition}
\newtheorem{lemma}[theorem]{Lemma}

\newtheorem{remark}[theorem]{Remark}

\numberwithin{equation}{section}
\numberwithin{theorem}{section}

\begin{document}

\title{Symmetric Ground States Solutions of m-Coupled nonlinear Schr\"odinger equations}
\author{Hichem Hajaiej}

\address{Justus-Liebig-Universit\"at Giessen\\Mathematisches Institut\\Arnd Str 2, 35392 Giessen\\Germany}\email{hichem.hajaiej@gmail.com}

\begin{abstract}
We prove the existence of radial and radially decreasing ground states of an m-coupled nonlinear Schr\"odinger equation with a general nonlinearity. 
\end{abstract}

\maketitle
\pagestyle{headings}
\section{Introduction}
\thispagestyle{empty}
The following Cauchy problem of an m-coupled \nlss: 
\begin{equation}\label{eq1.1}
\begin{cases}
i\partial_t\Phi_1+\Delta\Phi_1+g_1\left(|x|, |\Phi_1|^2, \ldots, |\Phi_m|^2\right)\Phi_1&= 0,\\
\qquad\qquad \vdots &\\
i\partial_t\Phi_m+\Delta\Phi_m+g_m\left(|x|, |\Phi_1|^2, \ldots, |\Phi_m|^2\right)\Phi_m&= 0,\\
\hfill\Phi_i(0,x)&=\Phi_i^0(x)\, \quad\mbox{for}\, 1\leq i\leq m. 
\end{cases}
\end{equation}
For $1\leq i\leq m: \Phi_i^0: \R^N\rightarrow\C$ and $g_i:\R_+^*\times\R_+^m\rightarrow\R$, 
$\Phi_i: \R_+\times\R^N\rightarrow\C$, has numerous applications in physical problems.  It appears in the study of spatial solitons in nonlinear waveguides~\cite{c22}, the theory of Bose-Einstein condensates~\cite{c6}, interactions of m-wave packets~\cite{c5}, optical pulse propagation in birefringent fibers~\cite{c17,c18}, wavelength division multiplexed optical systems.  Physically, the solution $\Phi_i$ is the $i$th component of the beam in Kerr-like photorefractive media~\cite{c1}.  In the most relevant cases, it is possible to write~(\ref{eq1.1}) in a vectorial form as follows: 
\begin{equation}\label{eq1.2}
\begin{cases}
i\frac{\partial\Phi}{\partial t}=E'(\Phi)&\\
\Phi(0,x)=\Phi^0=(\Phi^0_1, \ldots, \Phi^0_m)&
\end{cases}
\end{equation}
where
\begin{equation}
\label{eq1.3}
E(\Phi)=\frac{1}{2}\|\nabla\Phi\|_2^2-\int{G\left(|x|, \Phi_1, \ldots, \Phi_m\right)}\, dx. 
\end{equation}
$G: (0, \infty)\times \R^m\rightarrow\R$ satisfies the following system: 
\begin{equation}\label{eq1.4}
\begin{cases}
\frac{\partial G}{\partial u_1}=g_1\left(|x|, u_1^2, \ldots, u_m^2\right) u_1,&\\
\qquad \vdots &\\
\frac{\partial G}{\partial u_m}=g_m\left(|x|, u_1^2, \ldots, u_m^2\right) u_m.& 
\end{cases}
\end{equation}

When $m=1$, $G$ can be easily given by the explicit expression: $G(r,s)=\half\int_0^{s^2}g(r,t)\, dt$. 

In the general case: 
\begin{eqnarray}
G(r, u_1, \ldots, u_m)&=&\half\int_0^{u_1^2}g_1(r,t,u_2^2,\ldots,u_m^2)\, dt+K_1(u_2,\ldots,u_m)\nonumber\\
&=&\half\int_0^{u_i^2}g_i(r,u_1^2,\ldots,t_i,\ldots,u_m^2)\, dt_i+K_i(u_1,\ldots,u_{i-1},u_{i+1},\ldots,u_m)\nonumber\\
&=&\quad\ldots\nonumber\\
&=&\half\int_0^{u_m^2}g_m(r,u_1^2,\ldots,t)\, dt+K_m(u_1,\ldots,u_{m-1}).\label{eq1.5}
\end{eqnarray}

A soliton or standing wave of~(\ref{eq1.1}) is a solution of the form: $\Phi(t,x)=\left(\Phi_1(t,x), \ldots, \Phi_m(t,x)\right)$, where for $1\leq j\leq m: \Phi_j(t,x)=u_j(x)e^{-i\lambda_jt}$, $\lambda_j$ are real numbers.  Therefore $\U=(u_1, \ldots, u_m)$ is a solution of the following $m\times m$ elliptic eigenvalue problem: 
\begin{equation}\label{eq1.6}
\begin{cases}
\Delta u_1+\lambda_1 u_1+g_1\left(|x|, u_1^2, \ldots, u_m^2\right)u_1=0,&\\
\qquad \vdots & \\
\Delta u_m+\lambda_m u_m+g_m\left(|x|, u_1^2, \ldots, u_m^2\right)u_m=0.&
\end{cases}
\end{equation}

Among all the standing waves, let us mention the ground states which correspond to the least energy solutions $\U$ of~(\ref{eq1.6}), defined by: 
\begin{equation}\label{eq1.7}
E(\U)=\half\sum\limits_{i=1}^{m}|\nabla u_i|_2^2-\int_{\R^N}G\left(|x|,u_1(x),\ldots,u_m(x)\right)\, dx
\end{equation}
under constraints 
\begin{equation}\label{eq1.8}
S_c=\left\{\U=(u_1, \ldots, u_m)\in \mathrm{H}^1(\R^N)\times\ldots \mathrm{H}^1(\R^N): \int_{\R^N}u_i^2=c_i\right\}
\end{equation}
where $c_i>0$ are $m$ prescribed numbers. 

Ground states are solutions of the minimization problem: 
\begin{equation}\label{eq1.9}
\mbox{For given }\hfill c_i>0, M_c=\inf_{\U\in S_c}E(\U).\hfill
\end{equation}

Profiles of stable electromagnetic waves traveling along a medium are given by~(\ref{eq1.9}).  Note that in~(\ref{eq1.7}), $|x|$ is the position relative to the optical axis, $G$ is related to the index of refraction of the medium.  In the most relevant cases, $G$ has jumps at interfaces between layers of different media (core and claddings).  Therefore, $G$ is not continuous with respect to the first variable in many practical cases. 

The existence of ground states has been investigated by many authors following different methods.  In~\cite{c2,c8,c9,c13,c14new,c15new,c19,c23,c24,c25,c26new} by numerical arguments; in~\cite{c3,c4,c14,c15,c16,c20}, the mathematical analysis using the variational approach has been pursued to prove the existence of ground states.  These works addressed the special case $m=2$ and 
\begin{equation}\label{eq1.10}
\begin{cases}
g_1(|x|, u_1^2, u_2^2)=\left(|u_1|^{2p-2}+\beta|u_1|^{p-2}|u_2|^p\right),&\\
g_2(|x|, u_1^2, u_2^2)=\left(|u_2|^{2p-2}+\beta|u_2|^{p-2}|u_1|^p\right).&
\end{cases}
\end{equation}

This is a very interesting case where we can easily determine $G$, indeed using~(\ref{eq1.5}) it is obvious that $G(r, s_1, s_2)=\frac{1}{2p}u_1^{2p}+\frac{\beta}{p}u_1^pu_2^p+K_1(u_2)=\frac{1}{2p}u_2^{2p}+\frac{\beta}{p}u_1^pu_2^p+K_2(u_1)$.  A straightforward computation implies: $G(r,s_1,s_2)=\frac{1}{2p}s_1^{2p}+\frac{1}{2p}s_2^{2p}+\frac{\beta}{p}s_1^ps_2^p$. 

In~\cite{c3,c16}, not only the existence of ground states has been established, for~(\ref{eq1.1}) with $g_i$ given by~(\ref{eq1.10}), but also the orbital stability has been discussed.  Of course, we are interested in the orbital stability of ground states of~(\ref{eq1.1}) with general non-linearities.  However, an inescapable step consists in the establishment of suitable assumptions of $g_i$ under which~(\ref{eq1.1}) admits a unique solution.  This is a very challenging open question under investigation. 

Following a self-contained approach, we establish the existence of radial and radially decreasing ground states [Theorem~\ref{th3.1}].  Our main assumptions are that $G$ satisfies a growth condition and it is a supermodular function, that is to say: 
\begin{equation}\label{eq1.11}
G(r,y+he_i+ke_j)+G(r,y)\geq G(r,y+he_i)+G(r,y+ke_j)
\end{equation}
\begin{equation}\label{eq1.12}
G(r_1,y+he_i)+G(r_0,y)\leq G(r_1,y)+G(r_0,y+he_i)
\end{equation}
for every $i\neq j, h, k>0; y=(y_1, \ldots, y_m)$ and $e_i$ denotes the $i$th standard basis vector in $\R^m, r>0$ and $0<r_0<r_1$. 

These inequalities are connected to the cooperativity of~(\ref{eq1.6}). When $\lambda_i\equiv 0$, W.C. Troy proved in~\cite{c26} the necessity of this hypothesis.  Contrary to previous works, we will not use minimization under the so-called Nehari Manifold; neither results involving the Palais-Smale condition.  Instead, we take advantage of some recent results of symmetrization inequalities.  More precisely, in~\cite{c7}, it has been proved that if $G$ satisfies~(\ref{eq1.11}) and~(\ref{eq1.12}), then: 
\begin{equation}\label{eq1.13}
\int_{\R^N}G\left(|x|, u_1(x), \ldots, u_m(x)\right)\, dx \leq \int_{\R^N}G\left(|x|, u^*_1(x), \ldots, u^*_m(x)\right)\, dx. 
\end{equation}
Here $u^*$ denotes the Schwarz symmetrization of a function $u$ vanishing at infinity.  
It is well known that the norm of the  gradient does not increase under Schwarz symmetrization in L$^2$. Moreover rearrangements preserve the L$^2$~norm:
\begin{equation*}
\int|\nabla u^*|^2\leq\int|\nabla u|^2
\end{equation*}
\begin{equation*}
\int u^2=\int{(u^*)}^2. 
\end{equation*}

Finally let us point out that, as mentioned in~\cite{c34}, in many valuable papers the study of~(\ref{eq1.9}) with $m=1$ relied on the fact that one could look for minima in the class of radial functions using rearrangement inequalities. 
The compact embedding of such functions in~L$^p$ enables us to conclude~\cite{c27,c28,c29,c30,c31,c32,c33}.  H. Brezis and E.H. Lieb~\cite{c34} concluded this remark saying ``It is not known whether the minimum action lies in the class of radial solutions for $m>1$ because rearrangement inequalities are not applicable.''  In this paper we build on a method enabling us to use such vectorial inequalities to solve~(\ref{eq1.9}).

Thanks to these inequalities, we first prove that: Given $c_1, \ldots, c_m>0$: 
\begin{enumerate}
 \item (\ref{eq1.9}) always admits a minimizing sequence $\U_n=(u_{n,1}, \ldots, u_{n,m})$ such that each component $u_{n,i}$ is radial and radially decreasing. 
\item Noticing that any minimizing sequence of~(\ref{eq1.9}) is bounded, we will prove that if  $\U_n=\U^*_n\rightharpoonup\U$ then 
\begin{equation*}
\lim\limits_{n\rightarrow +\infty}\int_{\R^N}G\left(|x|, u_{n,1}, \ldots, u_{n,m}\right)\, dx=\int_{\R^N}G\left(|x|, u_1(x), \ldots, u_m(x)\right)\, dx
\end{equation*}
which implies that $\U=(u_1, \ldots, u_m)$ is such that $E(\U)\leq M_c$. 
\item To conclude, it is sufficient to prove that $\U\in S_c$. 
\end{enumerate}

This paper contains four more sections.  In the next section, we introduce the notation and definitions.  In the third section, we state our main result and give a detailed proof.  The fourth part is dedicated to a variant of our approach.  The last section is dedicated to some challenging open problems. 

\section{Preliminaries and Notation}
\begin{itemize}
\item In the sequel, $m, N\in\N^*$. 
\item For $1\leq p < \infty$, $|\cdot|_p$ denotes the norm in L$^p(\R^N)$. 
\item If $V=(v_1, \ldots, v_m)$ with $v_i\in$L$^2(\R^N): \|V\|_2^2=|v_1|_2^2+\ldots+|v_m|_2^2$. 
\item If $V=(v_1, \ldots, v_m)$ with $v_i\in$H$^1(\R^N): \|\nabla V\|_2^2=|\nabla v_1|_2^2+\ldots+|\nabla v_m|_2^2$. 
\item[] $[$H$^1(\R^N)]^m=$H$^1(\R^N)\times\ldots\times$H$^1(\R^N)$. 
\item All statements about measurability refer to the Lebesgue measure, $\mu$, on $\R^N$ or $(0, \infty)$.  When no domain of integration is indicated, the integral extends over $\R^N$. 
\item M$(\R^N)$ is the set of measurable functions on $\R^N$. 
\item F$(\R^N)$ is the set of symmetrizable functions: 
\begin{equation*}
\left\{u\in \textrm{M}(\R^N): u\geq 0\mbox{ and }\mu\{x\in\R^N: u(x)>t\}<\infty \quad\forall t>0\right\}. 
\end{equation*}
\item For $u\in$F$(R^N)$, $u^*$ denotes the Schwarz symmetrization of $u$.  For more details, see~\cite{c7}. 
\item We say that $u$ is Schwarz symmetric if $u\equiv u^*$. 
\item For $V\in\mathrm{F}(\R^N)\times\ldots\times\mathrm{F}(\R^N)$, $V$ is Schwarz symmetric if each of its components has its property. 
\item For the convenience of the reader, let us recall some important symmetrization inequalities~\cite{c10}: 
\begin{eqnarray}
\forall u\in\mathrm{H}^1(\R^N): |\nabla u|_2^2&=&\Big|\nabla|u|\Big|^2_2\geq\Big|\nabla|u|^*\Big|^2_2\label{eq2.1}\\
\forall u\in\mathrm{L}^2(\R^N): |u|_2^2&=&|u^*|_2^2. \label{eq2.2}
\end{eqnarray}
\end{itemize}
\begin{definition}\label{def2.1}
A function $G: (0,\infty)\times\R^m\rightarrow\R$ is an m-Carath\'eodory function if
\begin{enumerate}
 \item $G(\cdot, s_1, \ldots, s_m): (0,\infty)\rightarrow\R$ is measurable on $(0,\infty)\setminus\Gamma$, where $\Gamma$ is a subset of $(0,\infty)$ having one dimensional measure zero, for all $s_1, \ldots, s_m\geq 0$, 
\item For all $1\leq n\leq m$, every $(m-1)$ tuple $s_i\geq 0$ and $r\in(0,\infty)\setminus\Gamma$, the function: 
\begin{eqnarray}
\R&\rightarrow&\R\nonumber\\
s_n&\mapsto&G(r, \ldots, s_n, \ldots)\nonumber
\end{eqnarray}
is continuous on $\R$. 
\end{enumerate}
\end{definition}

This definition establishes the standard context for handling the measurability of the composite functions $G\left(|x|, u_1(x), \ldots, u_m(x)\right), u_i\in M(\R^N)$. 
An important property of an m-Carath\'eodory function is that $x\mapsto G\left(|x|, u_1(x), \ldots, u_m(x)\right)$ is measurable on $\R^N$ for every $u_1,\ldots,u_m\in M(\R^N)$
\begin{itemize}
 \item For the convenience of the reader, let us also recall that for an m-Carath\'eodory function satisfying~(\ref{eq1.11}) and~(\ref{eq1.12}), we have~(\ref{eq1.13}); \cite{c7}. 
\item For $r>0: B(0,r)=\{x\in\R^N:|x|<r\}, |x|$ is the euclidean norm in $\R^N$, there is a constant $V_N$ such that $\mu(B(0,r))=V_N r^N$ for all $r>0$. 
\end{itemize}

\section{Main result}\begin{theorem}\label{th3.1}
Let $G: (0,\infty)\times\R^m\rightarrow\R$ be such that: 
\begin{itemize}
\item[(G0)] $G$ is an m-Carath\'eodory function such that 
\begin{equation*}
G(r, s_1, \ldots, s_m)\leq G(r, |s_1|, \ldots, |s_m|)
\end{equation*}
for every $r>0$ and $s_1, \ldots, s_m\in\R$, 
\item[(G1)] For all $r>0; s_1, \ldots, s_m\geq 0$, we have
\begin{equation*}
 0\leq G(r, s_1, \ldots, s_m)\leq K\left(|s|^2+\sum\limits_{i=1}^{m}s_i^{\ell_i+2}\right):\newline
s=(s_1,\ldots,s_m); K>0 \textrm{ and } 0<\ell_i<\frac{4}{N}, 
\end{equation*}
\item[(G2)] $G$ satisfies~(\ref{eq1.11}) and ~(\ref{eq1.12}), 
\item[(G3)] $\forall\epsilon>0,\exists R_0>0$ and $S_0>0$ such that
$G(r, s_1, \ldots, s_m)\leq\epsilon|s|^2$ 
for all $r>R_0$, $s_1, \ldots, s_m<S_0; s=(s_1, \ldots, s_m)$,  
\item[(G4)] $G(r, t_1s_1, \ldots, t_ms_m)\geq t^2_{\max} G(r, s_1, \ldots, s_m)$ for any $t_1, \ldots, t_m\geq 1; r>0; s_1, \ldots, s_m\geq 0$ where $t_{\max}=\max\limits_{1\leq i\leq m}t_i$. 
\end{itemize}
Suppose additionally that $M_c<0$, then: \newline
$\forall c_1, \ldots, c_m>0$ there exist $V_c=\left(v_1^{c_1}, \ldots, v_m^{c_m}\right)$ such that $V_c\in S_c$ and $E(V_c)=M_c$. 
\end{theorem}

The proof of the result is divided in three parts: (step 1 $\rightarrow$ step 3): 
\begin{lemma}\label{lem3.1}
Suppose that $G$ satisfies (G0) and (G1), then all the minimizing sequences of~(\ref{eq1.9}) are bounded in $[\mathrm{H}^1(\R^N)]^m$. 
\end{lemma}

Proof: Let $\U=(u_1, \ldots, u_m)\in S_c$, (G0) and (G1) imply that 
\begin{equation*}
\int G(|x|,\U(x))\, dx\leq Kc+K\int\sum\limits_{i=1}^m|u_i(x)|^{\ell_i+2}\, dx.
\end{equation*}

For $1\leq i\leq m$, the Gagliardo-Nirenberg inequality tells us that: 
\begin{equation*}
|u_i|_{\ell_i+2}\leq C|u_i|_2^{1-\sigma_i}\cdot|\nabla u_i|_2^{\sigma_i}; \sigma_i=\frac{N}{2}\frac{\ell_i}{\ell_i+2}. 
\end{equation*}

Now let $\epsilon>0, p_i=\frac{4}{N\ell_i}, q_i$ is such that $\frac{1}{p_i}+\frac{1}{q_i}=1$. 
Applying Young's inequality, we obtain: 
\begin{equation*}
|u_i|_{\ell_i+2}\leq
\left\{
\frac{C^{\ell_i+2}}{\epsilon}|u_i|_2^{(1-\sigma_i)(\ell_i+2)}
\right\}^{q_i}
\frac{1}{q_i}+\frac{N\ell_i}{4}
\left\{
\epsilon^{\frac{4}{N\ell_i}} 
|\nabla u_i|^2_2\right\}.
\end{equation*}

Consequently: 
\begin{equation*}
E(\U)\geq \left\{\half-Km\sum_{i=1}^m\frac{N\ell_i}{4}\epsilon^{\frac{4}{N\ell_i}}\right\}\|\nabla\U\|_2^2-Kc-\sum\limits_{i=1}^m\frac{1}{q_i}C^{\ell_i+2}c^{\frac{(1-\sigma_i)(\ell_i+2)}{2}}. 
\end{equation*}

Taking $\epsilon$ such that $\half-Km\sum_{i=1}^m\frac{N\ell_i}{4}\epsilon^{\frac{4}{N\ell_i}}\geq0$, we prove that $E$ is bounded from below.  To show that any minimizing sequence of~(\ref{eq1.9}) is bounded in $[\mathrm{H}^1(\R^N)]^m$, it is enough to take the latter inequality with the strict sign. 
\begin{remark}\label{r3.1}
~\newline
\begin{itemize}
 \item The lemma remains true if we replace~(G1) by the more general growth condition: 
\begin{equation*}
G(r,s_1,\ldots,s_m)\leq K\left(|s|^2+\sum\limits_{k=0}^{\alpha}\left(\xi_{1,k}s_1+\ldots+\xi_{m,k}s_m\right)^{\ell_k+2}\right), 
\end{equation*}
for all $r>0$ and $s_1,\ldots,s_m\geq 0$, where $K$ is a positive constant, $\alpha\in\N^*$ and for $0\leq k\leq\alpha, 0<\ell_k<\frac{4}{N}$.  For $0\leq k\leq\alpha, 1\leq j\leq m: \xi_{j,k}$ can take arbitrarily the value 0 or 1. 
\item The growth condition stated in our lemma is optimal, in the sense that if $\ell>\frac{4}{N}$, we can prove that $M_c=-\infty$. 
\end{itemize}
\end{remark}

Under the hypotheses of Theorem~\ref{th3.1}, we will first prove that:\newline
{\bfseries Step 1:} 
\begin{equation}\label{eq3.1}
\mbox{For any }\U=(u_1, \ldots, u_m)\in\left[\mathrm{H}^1(\R^N)\right]^m: E(\U)\geq E(\U^*). 
\end{equation}

This inequality enables us to assert that for any m-tuple $c_1, \ldots, c_m>0$,~(\ref{eq1.9}) always admits a Schwarz symmetric minimizing sequence.  For such minimizing sequence, we have the following compactness property: \newline 
{\bfseries Step 2:} If $\U_n=U_n^*\rightharpoonup\U$ in $\left[\mathrm{H}^1(\R^N)\right]^m: E[\U]\leq\lim\inf E(\U_n)$. 

Finally we will show that this $\U$ belongs to the constraint when $M_c<0$. 

{\bfseries Step 1:} 
\begin{lemma}\label{lem3.2}
Suppose that $G$ satisfies (G0), (G1) and (G2).  If $(\U_n)$ is a minimizing sequence of~(\ref{eq1.9}), $\left(|\U_n|^*\right)$ also has this property. 
\end{lemma}

{\bfseries Proof:} Let $\U=(u_1, \ldots, u_m)\in\left[\mathrm{H}^1(\R^N)\right]^m$.  First note that for any $u_i\in\mathrm{H}^1(\R^N)$ and $|\nabla u_i|_2=\Big|\nabla|u_i|\Big|_2$, thus using (G0); $E(|\U|)=E(|u_1|, \ldots, |u_m|)\leq E(u_1, \ldots, u_m)$. 

To achieve the proof, it is sufficient to show that for any $V=(v_1, \ldots, v_m)$ with $v_i\geq 0$, $E(v_1^*, \ldots, v_m^*)\leq E(v_1, \ldots, v_m)$, which follows immediately from~(\ref{eq2.1}) and~(\ref{eq1.13}).  Note finally that by~(\ref{eq2.2}): if $\int v_i^2=c_i$ then $\int (v_i^*)^2=c_i$, this completes the proof. 

{\bfseries From now on:} 
\begin{equation}\label{eq3.2}
\U_n=(u_{n,1}, \ldots, u_{n,m})\textrm{ is a minimizing sequence of~(\ref{eq1.9}), which is Schwarz symmetric. }
\end{equation}
By Lemma~\ref{lem3.1}, it is bounded in $[\mathrm{H}^1(\R^N)]^m$.  We know that (up to a subsequence) there exists $\U=(u_1, \ldots, u_m)$ such that 
\begin{equation}\label{eq3.3}
u_{n,j}\rightharpoonup u_j \quad \forall 1\leq j\leq m. 
\end{equation}

{\bfseries Step 2:}
\begin{lemma}\label{lem3.3}
Let $G$ be a function satisfying (G0), (G1) and (G3). ($\U_n$) be a minimizing sequence satisfying~(\ref{eq3.2}) and~(\ref{eq3.3}) then $E(\U)\leq\lim\inf E(\U_n)$. 
\end{lemma}

{\bfseries Proof:} $\forall 1\leq i\leq m$, we know that $|\nabla u_i|_2^2\leq|\nabla u_{n,i}|_2^2$.  Let us prove that 
\begin{equation*}
\lim\limits_{n\rightarrow +\infty}\int G(|x|,u_{n,1}(x), \ldots, u_{n,m}(x))\, dx=\int G(|x|,u_{1}(x), \ldots, u_{m}(x))\, dx.
\end{equation*}

Let $R>0$, we first show that 
\begin{equation*}
\lim\limits_{n\rightarrow +\infty}\int_{|x|\leq R} G(|x|,u_{n,1}(x), \ldots, u_{n,m}(x))\, dx=\int_{|x|\leq R} G(|x|,u_{1}(x), \ldots, u_{m}(x))\, dx.
\end{equation*}

For $1\leq i\leq m$, $(u_{n,i})$ converges weakly to $u_i$ in H$^1(\R^N)$, it then converges to $u_i$ in L$^{\ell_i+2}$ $(|x|\leq R)$.  Therefore, up to a subsequence (which we also denote by $u_{n,i}$), $u_{n,i}\rightarrow u_i$ for almost every $|x|\leq R$, $|u_{n,i}|<h_i$ where $h_i\in\mathrm{L}^{\ell_i+2} (|x|\leq R)$. 

Now using (G1): 
\begin{equation*}
G(|x|,u_{n,1}(x), \ldots, u_{n,m}(x)) \leq  K\left(\sum\limits_{i=1}^mh^2_i(x)+\sum\limits_{i=1}^{m}h_i^{\ell_i+2}(x)\right).
\end{equation*}

All functions involved in this sum are in L$^1(|x|\leq R)$.  By the dominated convergence theorem, it follows that 
\begin{equation}\label{eq3.4}
\lim\limits_{n\rightarrow +\infty}\int_{|x|\leq R} G(|x|,u_{n,1}(x), \ldots, u_{n,m}(x))\, dx=\int_{|x|\leq R} G(|x|,u_{1}(x), \ldots, u_{m}(x))\, dx.
\end{equation}

Now fix $n\in\N$ and $1\leq i\leq n$.  Since $u_{n,i}$ is Schwarz symmetric: 
\begin{equation*}
V_N|x|^Nu_{n,i}^2(x)\leq\int_{|y|\leq|x|}u^2_{n,i}(y)\, dy\leq c_i.
\end{equation*}

Consequently $u_{n,i}(x)\leq\frac{c_i^{1/2}}{V_N^{1/2}|x|^{N/2}}\leq\frac{c_i^{1/2}}{V_N^{1/2}R^{N/2}}$ for all $|x|>R$. 

Let $\epsilon>0$, choose $R$ large enough, (G3) implies that 
\begin{equation*}
\int_{|x|>R}G(|x|,u_{n,1}(x), \ldots, u_{n,m}(x))\, dx\leq\epsilon\sum\limits_{i=1}^m\int_{|x|>R}u_{n,i}^2(x)\, dx\leq \epsilon c,
\end{equation*}
where $c=\sum\limits_{i=1}^mc_i$. 

Proving that: 
\begin{equation}\label{eq3.5}
\lim\limits_{R\rightarrow\infty}\lim\limits_{n\rightarrow\infty}\int\limits_{|x|>R}G(|x|,u_{n,1}(x), \ldots, u_{n,m}(x))\, dx=0.
\end{equation}

The two properties we need to prove~(\ref{eq3.5}) are: $\int u_{n,i}^2(x)\leq c_i$ and $(u_{n,i})$ is Schwarz symmetric $\forall 1\leq i\leq m$. 

Clearly $\int u_i^2\leq c_i$.  The second property is inherited by $u_i$ almost everywhere.  Indeed for $R>0$, there exists $n_k(R)$ such that $(u_{n_k,i})$ converges to $u_i$ almost everywhere and we obtain: 
\begin{equation*}
\lim\limits_{R\rightarrow\infty}\int\limits_{|x|>R}G(|x|,u_{1}(x), \ldots, u_{m}(x))\, dx=0. 
\end{equation*}
Consequently 
\begin{equation*}
\lim\limits_{n\rightarrow\infty}\int G(|x|,u_{n,1}(x), \ldots, u_{n,m}(x))\, dx=\int  G(|x|,u_{1}(x), \ldots, u_{m}(x))\, dx. 
\end{equation*}

Thanks to our lemmas, we know that $E(\U)\leq M_c$; ($\U=(u_1, \ldots, u_m)$ is given by~(\ref{eq3.3})): 
\begin{equation}\label{eq3.6}
|u_i|^2_2\leq c_i \quad \forall 1\leq i\leq m. 
\end{equation}

{\bfseries Step 3:} To conclude that the infinum is achieved, we have to prove that $\U\in S_c$.  Suppose that $M_c<0$, set $t_i=\frac{c_i^{1/2}}{|u_i|_2}$, by~(\ref{eq3.6}): 
\begin{equation}\label{eq3.7}
t_i\geq 1\mbox{ and }(t_1u_1, \ldots, t_mu_m)\in S_c \qquad t_{\max}=\max\limits_{1\leq i\leq m}t_i\geq 1.  
\end{equation}
\begin{equation*}
E\left(t_1u_1, \ldots, t_mu_m\right)=\half\sum\limits_{i=1}^m|t_i\nabla u_i|^2_2-\int G(|x|, t_1u_1(x), \ldots, t_mu_m(x))\, dx. 
\end{equation*}

By (G4): 
\begin{eqnarray*}
E(t_1u_1, \ldots t_mu_m)&\leq& t_{\max}^2E(u_1, \ldots, u_m). \\
M_c\leq E(t_1u_1, \ldots t_mu_m)&\leq& t_{\max}^2E(u_1, \ldots, u_m) \leq t_{\max}^2M_c,  
\end{eqnarray*}
since $t_i\geq 1$ by Lemma~\ref{lem3.3}, it follows that $M_c\leq t^2_{\max}M_c \Rightarrow t^2_{\max}\leq 1$, hence $t_i=1$ for any $1\leq i\leq m$.  This ends the proof of Theorem~\ref{th3.1}. 

{\bf On the hypothesis $\mathbf{M_c<0}$:} \newline
Inspired by~\cite{c21} and closely following the approach therein, we prove that if $G$ satisfies: 
\begin{itemize}
 \item[(G5)] There exist $R_1>0, S_1>0$.  For any $1\leq i\leq m$, there exist $A_i>0, t_i\in [0,2)$ and $0\leq\sigma_i\leq\frac{2(2-t_i)}{N}$ such that
\begin{equation*}
G(r,s_1,\ldots,s_m)\geq \sum\limits_{i=1}^mA_ir^{-t_i}s_i^{\sigma_i+2}\mbox{ for all }r>R_1, 0<s<s_1
\end{equation*}
then $M_c<0$. 
\end{itemize}

Set $d(N)=\int e^{-2|y|^2}\, dy, D(N)=\frac{4}{d^2(N)}\int|y|^2e^{-2|y|^2}\, dy$.  For $\alpha\in (0,1]$, we set $w_\alpha: \R^N\rightarrow\R$ defined by $w_\alpha(x)=\frac{\alpha^{N/4}e^{-\alpha|x|^2}}{d(N)}$.  A straightforward computation shows that $|w_\alpha|_2=1$ and $|\nabla w_\alpha|_2^2=\alpha D(N)$. 

On the other hand, there exists $B>R_1$ such that for any $|x|>B$, $w_\alpha(x)\leq S_1$. 
\begin{equation*}
\int G(|x|,w_\alpha(x),\ldots,w_\alpha(x))\geq \int_{|x|\geq B}\sum\limits_{i=1}^{m}\frac{A_i}{[d(N)]^{\sigma_i+2}}|x|^{-t_i}e^{-\alpha(\sigma_i+2)|x|^2}\alpha^{\frac{N}{4}(\sigma_i+2)}\, dx. 
\end{equation*}
By the change of variable $y=\alpha^\half x$, we obtain: 
\begin{eqnarray*}
&=&\sum\limits_{i=1}^{m}\frac{A_i}{[d(N)]^{\sigma_i+2}}\alpha^{\frac{N\sigma_i}{4}+\frac{t_i}{2}}\int_{|y|\geq B\alpha^\half}|y|^{-t_i}e^{-(\sigma_i+2)|y|^2}\, dy
\\
&\geq&\sum\limits_{i=1}^{m}\frac{A_i}{[d(N)]^{\sigma_i+2}}\alpha^{\frac{N\sigma_i}{4}+\frac{t_i}{2}}\int_{|y|\geq B}|y|^{-t_i}e^{-(\sigma_i+2)|y|^2}\, dy
\end{eqnarray*}

Set $I_i=\int_{|y|\geq B}|y|^{-t_i}e^{-(\sigma_i+2)|y|^2}\, dy$, it follows that: 
\begin{equation*}
E(w_\alpha, \ldots, w_\alpha)\leq \alpha\left\{mD(N)-\sum\limits_{i=1}^m\frac{A_i}{[d(N)]^{\sigma_i+2}}I_i\alpha^{\frac{N\sigma_i}{4}+\frac{t_i}{2}-1} \right\}. 
\end{equation*}

The fact that $\sigma_i<2(2-t_i)/N$ enables us to conclude that $E(w_\alpha, \ldots, w_\alpha)<0$ for $\alpha$ sufficiently small.  Taking $u_i=\frac{c_i^{1/2} w_\alpha}{|w_\alpha|_2}$, we can easily see that $E(u_1, \ldots, u_m)<0$ with $(u_1, \ldots, u_m)\in S_c$, thus $M_c<0$. 

\section{Variant of our result}
Our approach also applies to the following variational problem: 
\begin{equation*}
\tilde{M_c}=\inf_{\U\in\tilde{S_c}}\tilde{E}(\U),\mbox{ for } \U=(u_1, \ldots, u_m)\in\left[\mathrm{H}^1(\R^N)\right]^m, 
\end{equation*}
\begin{equation*}
\tilde{E}(\U)=\half\sum\limits_{i=1}^m|\nabla u_i|_2^2-\half\int p(|x|)\sum\limits_{i=1}^mu_i^2(x)-\int G(|x|,u_1(x),\ldots,u_m(x)). 
\end{equation*}

For a prescribed $c>0$: $\tilde{S_c}=\left\{\U=(u_1, \ldots, u_m): \|\U\|^2_2=c\right\}$.  Then we have the following result: 
\begin{theorem}\label{th4.1}
Suppose that $p:(0,\infty)\rightarrow\R$ satisfies
\begin{itemize}
\item[(P1)] $p$ is non-negative, non-increasing and $\lim\limits_{r\rightarrow\infty}p(r)=0$; 
\item[(P2)]\begin{itemize}
\item If $N=1, 2$, there exists $a\in (0,1]$ such that $p(a)>0$; 
\item If $N\geq 3$, there exists $R>0$ such that $p(r)>\frac{j^2_{N/2-1,1}}{R^2}$ where $j^2_{N/2-1,1}$ is the first zero of the Bessel function~$J_{N/2-1}$. 
\end{itemize}
\end{itemize}
Suppose that $G$ satisfies $(G0)\rightarrow (G4)$ in which each~$t_i$ is replaced by~$t$, then, for any $c>0$, there exists $\U_c=(u_c^1, \ldots, u_c^m)$ Schwarz symmetric such that~$\tilde{E}(\U_c)=\tilde{M_c}$. 
\end{theorem}

{\bf Proof:} Following the same approach as in the previous Theorem, step~1, step~2 and step~3 can be proven under minor modifications.  Therefore we are done if $\tilde{M_c}<0$.  Since~$G$ is non-negative, it is sufficient to prove that we can construct $v\in\mathrm{H}^1(\R^N)$ such that
\begin{equation}\label{eq4.1}
\half|\nabla v|_2^2-\half\int p(|x|)v^2<0.
\end{equation}

For the convenience of the reader, we will mention all the details.  These test functions were constructed in~\cite{c11} and used in~\cite{c12}. 

\begin{itemize}
\item Case $N=1$: \newline
Take $w(x)=e^{-|x|}$, $\alpha\in (0,1], 0<d\leq a$ and $w_\alpha(x)=w(\alpha x)$ 
\begin{equation}\label{eq4.2}
\half\int|\nabla w_\alpha|^2-p(|x|)w^2_\alpha(x)\, dx=\half\int\alpha^2|\nabla w(\alpha x)|^2-p(|x|)w^2(\alpha x)\, dx. 
\end{equation}
By the change of variables $y=\alpha x$, we obtain: 
\begin{eqnarray*}
(\ref{eq4.2}) &\leq& \frac{1}{2\alpha}\left\{\alpha^2|\nabla w|_2^2-\int p\left(\frac{|y|}{\alpha}\right)w^2(y)\, dy\right\}\leq
\frac{1}{2\alpha}\left\{\alpha^2|\nabla w|_2^2-w^2(d)\int_{|y|\leq d} p\left(\frac{|y|}{\alpha}\right)\, dy\right\}\\
(\ref{eq4.2}) &\leq& \frac{\alpha}{2}\left\{|\nabla w|_2^2-\frac{w^2(d)p(d)2d}{\alpha}\right\}. 
\end{eqnarray*}

In the last inequality, we have used the change of variables $z=\frac{y}{\alpha}$, then used the monotonicity of $p$.  

Therefore for $\alpha$ small enough,~(\ref{eq4.2})$<0$. Now for $c>0$ and $\alpha$ small enough take: $v_i=\frac{c^{1/2}w_\alpha}{m^{1/2}|w_\alpha|_2}$, then $\half\int|\nabla v_i|_2^2-\half\int p(|x|)v_i^2<0$, $v=(v_1, \ldots, v_m)\in\tilde{S_c}$ and $\tilde{E}(v_1, \ldots, v_m)<0$.  
\item Case $N=2$: \newline 
Let $
u(x)=
\begin{cases}
\left(\log{\frac{1}{|x|}}\right)^{1/3}&\mbox{if }|x|<1, \\
0&\mbox{otherwise.}
\end{cases}
$

$u\in\mathrm{H}^1(\R^2)$ but it is an unbounded function because of its singularity in~$0$.  Let $K=\left(\int_{|x|\leq 1}p(|x|)\, dx\right)^{-1}$, there exists $d\in\R^2$ such that 
\begin{equation}\label{eq4.3}
u^2(d)>K|\nabla u|_2^2. 
\end{equation}

Set $w_d(x)=u(|d|x), w_d\in\mathrm{H}^1(\R^2)$ and: 
\begin{eqnarray*}
&&\half|\nabla w_d|^2_2-\half\int p(|x|)w_d^2(x)\, dx\leq\half\int|d|^2\Big|\nabla u(|d|x)\Big|^2-p(|x|)u^2(|d|x)\, dx\\
&\leq&\half\int|\nabla u(y)|^2-\frac{1}{|d|^2}p\left(\frac{|y|}{|d|}\right)u^2(y)\, dy\leq
\half\int|\nabla u(y)|^2-\frac{1}{2|d|^2}\int_{|y|\leq d}p\left(\frac{|y|}{|d|}\right)u^2(y)\, dy\\
&\leq&\half\left\{|\nabla u|^2_2-u^2(d)\int_{|z|\leq 1}p(|z|)\, dz\right\}<0\mbox{ by (\ref{eq4.3})}. 
\end{eqnarray*}

The proof goes as previously setting~$v_i=\frac{c^{1/2}w_d}{m^{1/2}|w_d|_2}$ for $1\leq i\leq n$. 
\item Case $N\geq 3$: Let $x\in\B(0,1)$, set $\varphi_1(x)=|x|^{-\left(\frac{N}{2}-1\right)}J_{N/2-1}\left(j_{N/2-1,1}|x|\right)$.  It is easy to check that $\varphi_0\in\mathrm{H}_0^1(|x|<1)$ and $-\Delta\varphi_1=j^2_{N/2-1,1}\varphi_1$.  For $R$ given by~(P2), set $\varphi_R(x)=\varphi_1\left(\frac{x}{R}\right)$ then $\varphi_R\in\mathrm{H}_0^1(|x|<R)$ and $-\nabla\varphi_R=\frac{j^2_{N/2-1,1}}{R^2}\varphi_R$. 

Now set 
$w_R=
\begin{cases}
\varphi_R&\mbox{if }|x|<R\\
0&\mbox{otherwise.}
\end{cases}$

$w_R\in\mathrm{H}^1(\R^N)$ and $\half\int|\nabla w_R|-\half\int p(|x|)w^2_R(x)\, dx\leq\half\int_{|x|\leq R}\left\{\frac{j^2_{N/2-1,1}}{R^2}-p(|x|)\right\}w_R^2(x)\, dx<0$ by~(P2). 

We conclude in the same way as in the previous cases. 
\end{itemize}
\begin{remark}
Theorem~\ref{th4.1} holds true when~(P2) is replaced by~(G5). 
\end{remark}

\vspace{4ex}
{\bfseries Examples of functions $\mathbf{G}$ satisfying $\mathbf{(G0)\rightarrow(G5)}$:} 

Let $m=2, k\in\N^*$:\newline
(R)\qquad\qquad\qquad\qquad\qquad\qquad$G(r,s)=b(r)|s|^2+a(r)\sum\limits_{j=1}^k|s_1|^{\ell_{1,j}+1}|s_2|^{\ell_{2,j}+1}$

\begin{itemize}
 \item[(R1)] $\ell_{1,j}$ and $\ell_{2,j}>1$ with $\ell_{1,j}+\ell_{2,j}<\frac{4}{N}$ for $1\leq j\leq k$. 
 \item[(R2)] $a(r)$ is a non-negative, non-increasing function bounded from above and below by two positive constants. 
 \item[(R3)] $b(r)$ is a non-negative, non-increasing bounded function tending to zero as $r$ goes to infinity. 
\end{itemize}
Then $G$ satisfies $(G0)\rightarrow(G5)$. 

{\bfseries Remarks:} 
\begin{itemize}
 \item For $m>2$, functions $G$ satisfying $(G0)$ to $(G5)$ are given in a similar way as~(R) with a sum involving products of all $s_i, 1\leq i\leq m$.  This ensures~$(G4)$. 
\item Note that in~(R), $|s|^2$ can be replaced by~$|s|^{\sigma+2}$ with $0<\sigma<\frac{4}{N}$.  In this case~$b(r)$ can be taken as a positive constant: (R')
\item Finally when one deals with functions $G$ that are not necessarily sums of products involving all~$s_i$ with~$1\leq i\leq m$, we should apply Theorem~\ref{th4.1}, from which we can easily see that~(\ref{eq1.10}) is a particular case of this result.  More precisely, take~$a\equiv \frac{\beta}{p}$, $b=\frac{1}{2p}$, $\ell_1=\ell_2=\frac{\sigma}{2}=p-1$ with $1<p<\frac{2}{N}$ in~(R'). 
\end{itemize}

\section{Concluding remarks}
In this paper, we have determined suitable assumptions of the operator $G$, involved in the m-coupled nonlinear Schr\"odinger equations such that~(\ref{eq1.1}) admits a radial and radially decreasing ground state with respect to each component.  Moreover, if~(\ref{eq1.11}) and~(\ref{eq1.12}) hold true with strict inequality~\cite[Theorem 2]{c13}, it follows that $E(\U^*)<E(\U)$ for any~$\U\in[H^1(\R^N)]^m$.  Consequently all the ground states of~(\ref{eq1.1}) are Schwarz symmetric.  A challenging question is the establishment of the uniqueness of these least energy solutions.  Until now, we are not aware of any result in this direction when $N>1$ and $m>1$.  Another very interesting question is the study of the orbital stability of these standing waves.  We expect that for $\ell_i<4/N$, the ground states are stable.  A crucial step to establish such a result is to prove the uniqueness of the solutions of~(\ref{eq1.1}).  
For more general nonlinearities $g_i$, this open problem, under investigation, seems to be extremely complicated. 

\section*{Acknowledgment}
The author is extremely grateful to Dr. Yvan Pointurier for his precious help. 
The author is also grateful to the referees, Stefan Le Coz and Louis Jeanjean  for their  valuable comments.

\end{document}